\documentclass{amsart}
\usepackage{a4wide}
\usepackage{ascmac}
\usepackage{amsmath}
\usepackage{amscd}
\usepackage{amssymb}
\usepackage{hyperref}
\usepackage{graphicx,color}
\usepackage{pict2e}
\usepackage{here}

\theoremstyle{definition}
\newtheorem{definition}{Definition}
\newtheorem{example}{Example}
\theoremstyle{plain}
\newtheorem{claim}{Claim}

\newtheorem{theorem}{Theorem}
\newtheorem{corollary}{Corollary}
\newtheorem{proposition}{Proposition}
\newtheorem{fact}{Fact}

\theoremstyle{remark}
\newtheorem{remark}{Remark}

\newcommand{\rii}{\operatorname{RI\!I}}
\newcommand{\Gal}{ 
\begin{picture}(60,30)(0,0)
\put(12,12.5){\circle{20}}
\put(48,12.5){\circle{20}}
\put(20,7.5){\vector(2,1){20}}
\put(40,7.5){\vector(-2,1){20}}
\end{picture}
}
\newcommand{\GalS}{ 
\begin{picture}(60,30)(0,0)
\put(12,12.5){\circle{20}}
\put(48,12.5){\circle{20}}
\put(20,17.5){\vector(2,-1){20}}    
\put(40,17.5){\vector(-2,-1){20}} 
\end{picture}
}
\newcommand{\LK}{
\begin{picture}(30,10)
\put(5,3){\circle{10}}
\put(10,3){\vector(1,0){11}}
\put(26,3){\circle{10}}
\end{picture}
}
\newcommand{\Ga}{ 
\begin{picture}(60,30)(0,0)
\put(12,12.5){\circle{20}}
\put(48,12.5){\circle{20}}
\put(20,17.5){\vector(1,0){20}}    
\put(40,7.5){\vector(-1,0){20}}      
\end{picture}
}
\newcommand{\GaS}{
\begin{picture}(60,30)(0,0)
\put(12,12.5){\circle{20}}
\put(48,12.5){\circle{20}}
\put(40,17.5){\vector(-1,0){20}}   
\put(20,7.5){\vector(1,0){20}}   
\end{picture}
}
\newcommand{\Gb}{ 
\begin{picture}(60,30)(0,0)
\put(12,12.5){\circle{20}}
\put(48,12.5){\circle{20}}
\put(20,17.5){\vector(1,0){20}}    
\put(20,7.5){\vector(1,0){20}}   
\end{picture}
}
\newcommand{\GbS}{ 
\begin{picture}(60,30)(0,0)
\put(12,12.5){\circle{20}}
\put(48,12.5){\circle{20}}
\put(20,17.5){\vector(2,-1){20}}    
\put(20,7.5){\vector(2,1){20}}   
\end{picture}
}
\newcommand{\Gbl}{ 
\begin{picture}(60,30)(0,0)
\put(12,12.5){\circle{20}}
\put(48,12.5){\circle{20}}
\put(40,17.5){\vector(-1,0){20}}   
\put(40,7.5){\vector(-1,0){20}}      
\end{picture}
}
\newcommand{\GblS}{ 
\begin{picture}(60,30)(0,0)
\put(12,12.5){\circle{20}}
\put(48,12.5){\circle{20}}
\put(40,17.5){\vector(-2,-1){20}}   
\put(40,7.5){\vector(-2,1){20}}      
\end{picture}
}

\newcommand{\OR}{ 
\begin{picture}(50,40)(0,0)
\qbezier(5,27.5)(12.5,35)(20,27.5)  
\qbezier(5,12.5)(12.5,5)(20,12.5)    
\qbezier(40,2.5)(30,20)(40,37.5)     
\put(12.5,31.5){\line(1,0){25}}          
\put(12.5,8.5){\line(1,0){25}}            
\put(17,10){\line(0,1){20}}                
\end{picture}
}
\newcommand{\ORS}{ 
\begin{picture}(50,40)(0,0)
\qbezier(5,27.5)(12.5,35)(20,27.5)  
\qbezier(5,12.5)(12.5,5)(20,12.5)    
\qbezier(40,2.5)(30,20)(40,37.5)     
\qbezier(15,31)(17.5,34)(20,31.5)   
\qbezier(15,9)(17.5,6)(20,8.5)         
\put(20,31.5){\line(4,-5){16}}          
\put(20,8.5){\line(4,5){16}}            
\put(9,10){\line(0,1){20}}                
\end{picture}
}

\newcommand{\ORab}{ 
\begin{picture}(50,40)(0,0)
\qbezier(5,27.5)(12.5,35)(20,27.5)  
\qbezier(5,12.5)(12.5,5)(20,12.5)    
\qbezier(40,2.5)(30,20)(40,37.5)     
\put(12.5,31.5){\line(1,0){25}}          
\put(17,10){\line(0,1){20}}                
\end{picture}
}
\newcommand{\ORabS}{ 
\begin{picture}(50,40)(0,0)
\qbezier(5,27.5)(12.5,35)(20,27.5)  
\qbezier(5,12.5)(12.5,5)(20,12.5)    
\qbezier(40,2.5)(30,20)(40,37.5)     
\qbezier(15,31)(17.5,34)(20,31.5)   
\put(20,31.5){\line(4,-5){16}}          
\put(9,10){\line(0,1){20}}                
\end{picture}
}

\newcommand{\ORar}{ 
\begin{picture}(50,40)(0,0)
\qbezier(5,27.5)(12.5,35)(20,27.5)  
\qbezier(5,12.5)(12.5,5)(20,12.5)    
\qbezier(40,2.5)(30,20)(40,37.5)     
\put(12.5,8.5){\line(1,0){25}}            
\put(17,10){\line(0,1){20}}                
\end{picture}
}
\newcommand{\ORarS}{ 
\begin{picture}(50,40)(0,0)
\qbezier(5,27.5)(12.5,35)(20,27.5)  
\qbezier(5,12.5)(12.5,5)(20,12.5)    
\qbezier(40,2.5)(30,20)(40,37.5)     
\qbezier(15,9)(17.5,6)(20,8.5)         
\put(20,8.5){\line(4,5){16}}            
\put(9,10){\line(0,1){20}}                
\end{picture}
}

\newcommand{\ORbr}{ 
\begin{picture}(50,40)(0,0)
\qbezier(5,27.5)(12.5,35)(20,27.5)  
\qbezier(5,12.5)(12.5,5)(20,12.5)    
\qbezier(40,2.5)(30,20)(40,37.5)     
\put(12.5,31.5){\line(1,0){25}}          
\put(12.5,8.5){\line(1,0){25}}            
\end{picture}
}
\newcommand{\ORbrS}{ 
\begin{picture}(50,40)(0,0)
\qbezier(5,27.5)(12.5,35)(20,27.5)  
\qbezier(5,12.5)(12.5,5)(20,12.5)    
\qbezier(40,2.5)(30,20)(40,37.5)     
\qbezier(15,31)(17.5,34)(20,31.5)   
\qbezier(15,9)(17.5,6)(20,8.5)         
\put(20,31.5){\line(4,-5){16}}          
\put(20,8.5){\line(4,5){16}}            
\end{picture}
}

\newcommand{\Aright}{ 
\begin{picture}(50,40)(0,0)
\qbezier(5,27.5)(12.5,35)(20,27.5)  
\qbezier(5,12.5)(12.5,5)(20,12.5)    
\qbezier(40,2.5)(30,20)(40,37.5)     
\put(12.5,31.5){\vector(1,0){25}}          
\put(12.5,8.5){\vector(1,0){25}}            
\end{picture}
}
\newcommand{\ArightS}{ 
\begin{picture}(50,40)(0,0)
\qbezier(5,27.5)(12.5,35)(20,27.5)  
\qbezier(5,12.5)(12.5,5)(20,12.5)    
\qbezier(40,2.5)(30,20)(40,37.5)     
\qbezier(15,31)(17.5,34)(20,31.5)   
\qbezier(15,9)(17.5,6)(20,8.5)         
\put(20,31.5){\vector(4,-5){16}}          
\put(20,8.5){\vector(4,5){16}}            
\end{picture}
}

\newcommand{\Aleft}{ 
\begin{picture}(50,40)(0,0)
\qbezier(5,27.5)(12.5,35)(20,27.5)  
\qbezier(5,12.5)(12.5,5)(20,12.5)    
\qbezier(40,2.5)(30,20)(40,37.5)     
\put(37.5,31.5){\vector(-1,0){24}}          
\put(37.5,8.5){\vector(-1,0){24}}            
\end{picture}
}
\newcommand{\AleftS}{ 
\begin{picture}(50,40)(0,0)
\qbezier(5,27.5)(12.5,35)(20,27.5)  
\qbezier(5,12.5)(12.5,5)(20,12.5)    
\qbezier(40,2.5)(30,20)(40,37.5)     
\qbezier(15,31)(17.5,34)(20,31.5)   
\qbezier(15,9)(17.5,6)(20,8.5)         
\put(15,31){\vector(-4,-2){0.5}}
\put(15,9){\vector(-4,2){0.5}}
\put(20,31.5){\line(4,-5){16}}          
\put(20,8.5){\line(4,5){16}}            
\end{picture}
}
\newcommand{\Bright}{ 
\begin{picture}(50,40)(0,0)
\qbezier(5,27.5)(12.5,35)(20,27.5)
\qbezier(5,12.5)(12.5,5)(20,12.5)    
\qbezier(40,2.5)(30,20)(40,37.5)     
\put(12.5,31.5){\vector(1,0){25}}         
\put(37.5,8.5){\vector(-1,0){24}}            
\end{picture}
}
\newcommand{\BrightS}{ 
\begin{picture}(50,40)(0,0)
\qbezier(5,27.5)(12.5,35)(20,27.5)  
\qbezier(5,12.5)(12.5,5)(20,12.5)    
\qbezier(40,2.5)(30,20)(40,37.5)     
\qbezier(15,31)(17.5,34)(20,31.5)   
\qbezier(15,9)(17.5,6)(20,8.5)         
\put(15,9){\vector(-4,2){0.5}}
\put(20,31.5){\vector(4,-5){16}}
\put(20,8.5){\line(4,5){16}}            
\end{picture}
}

\newcommand{\GaMinusMinus}{ 
\begin{picture}(60,10)(0,12)
\put(12,22.5){\circle*{2}}
\put(12,12.5){\circle{20}}
\put(48,22.5){\circle*{2}}
\put(48,12.5){\circle{20}}
\put(27.5,20){\tiny$-$}
\put(20,17.5){\vector(1,0){20}}    
\put(27.5,10){\tiny$-$}
\put(40,7.5){\vector(-1,0){20}}      
\end{picture}
}
\newcommand{\GaSMinusMinus}{
\begin{picture}(60,10)(0,12)
\put(12,22.5){\circle*{2}}
\put(12,12.5){\circle{20}}
\put(48,22.5){\circle*{2}}
\put(48,12.5){\circle{20}}
\put(27.5,20){\tiny$-$}
\put(40,17.5){\vector(-1,0){20}}   
\put(27.5,10){\tiny$-$}
\put(20,7.5){\vector(1,0){20}}   
\end{picture}
}
\newcommand{\GaSMinusPlus}{
\begin{picture}(60,10)(0,12)
\put(12,22.5){\circle*{2}}
\put(12,12.5){\circle{20}}
\put(48,22.5){\circle*{2}}
\put(48,12.5){\circle{20}}
\put(27.5,20){\tiny$-$}
\put(40,17.5){\vector(-1,0){20}}   
\put(27.5,10){\tiny$+$}
\put(20,7.5){\vector(1,0){20}}   
\end{picture}
}

\newcommand{\GbPlusPlus}{ 
\begin{picture}(60,10)(0,12)
\put(12,22.5){\circle*{2}}
\put(12,12.5){\circle{20}}
\put(48,22.5){\circle*{2}}
\put(48,12.5){\circle{20}}
\put(27.5,20){\tiny$+$}
\put(20,17.5){\vector(1,0){20}}    
\put(27.5,10){\tiny$+$}
\put(20,7.5){\vector(1,0){20}}   
\end{picture}
}

\newcommand{\GblPlusPlus}{ 
\begin{picture}(60,10)(0,12)
\put(12,22.5){\circle*{2}}
\put(12,12.5){\circle{20}}
\put(48,22.5){\circle*{2}}
\put(48,12.5){\circle{20}}
\put(27.5,20){\tiny$+$}
\put(40,17.5){\vector(-1,0){20}}   
\put(27.5,10){\tiny$+$}
\put(40,7.5){\vector(-1,0){20}}      
\end{picture}
}
\newcommand{\GbMinusMinus}{ 
\begin{picture}(60,10)(0,12)
\put(12,22.5){\circle*{2}}
\put(12,12.5){\circle{20}}
\put(48,22.5){\circle*{2}}
\put(48,12.5){\circle{20}}
\put(27.5,20){\tiny$-$}
\put(20,17.5){\vector(1,0){20}}    
\put(27.5,10){\tiny$-$}
\put(20,7.5){\vector(1,0){20}}   
\end{picture}
}

\newcommand{\GblMinusMinus}{ 
\begin{picture}(60,10)(0,12)
\put(12,22.5){\circle*{2}}
\put(12,12.5){\circle{20}}
\put(48,22.5){\circle*{2}}
\put(48,12.5){\circle{20}}
\put(27.5,20){\tiny$-$}
\put(40,17.5){\vector(-1,0){20}}   
\put(27.5,10){\tiny$-$}
\put(40,7.5){\vector(-1,0){20}}      
\end{picture}
}

\newcommand{\GbMinusPlus}{ 
\begin{picture}(60,10)(0,12)
\put(12,22.5){\circle*{2}}
\put(12,12.5){\circle{20}}
\put(48,22.5){\circle*{2}}
\put(48,12.5){\circle{20}}
\put(27.5,20){\tiny$-$}
\put(20,17.5){\vector(1,0){20}}    
\put(27.5,10){\tiny$+$}
\put(20,7.5){\vector(1,0){20}}   
\end{picture}
}

\newcommand{\GblMinusPlus}{ 
\begin{picture}(60,10)(0,12)
\put(12,22.5){\circle*{2}}
\put(12,12.5){\circle{20}}
\put(48,22.5){\circle*{2}}
\put(48,12.5){\circle{20}}
\put(27.5,20){\tiny$-$}
\put(40,17.5){\vector(-1,0){20}}   
\put(27.5,10){\tiny$+$}
\put(40,7.5){\vector(-1,0){20}}      
\end{picture}
}

\title[Multiple linking number]{Multiple linking number}
\author{Kamolphat Intawong and Noboru Ito}
\begin{document}
\begin{abstract}
The linking number is the simplest link invariant given by Gauss; it is the first Gauss diagram formula expressed by one arrow among two circles.  Proceeding the next stage,  
we study the second Gauss diagram formula consisting of two arrows among two circles.  We call a function of this type the \emph{multiple linking number}.  There are two multiple linking numbers; one of them is ordinary Vassiliev invariant and the other function is surprisingly sensitive to  
 the necessity of the second  Reidemeister moves  
 though any one-component Gauss diagram formula cannot detect the necessity.
\end{abstract}
\date{\today}
\maketitle
\begin{sloppypar}
\section{Introduction}\label{sec:intro}
Gausss formulated the linking number $lk(K_1, K_2)$ of a link $K_1 \cup K_2$
\begin{align}
lk(K_1, K_2) &=\frac{1}{4\pi} \int_{S^1 \times S^1} \frac{(g(t) - f(s)) \times f'(s) \cdot  g'(t)}{|| f(s) - g(t) ||^3}  ds dt \label{eqGauss} \\
&= \int_{S^1} B(f(t)) \cdot  g'(t)  dt \quad, \nonumber
\end{align}
where $B(x)$ is a magnetic filed (the well-known Biot-Savart law).  
The Gauss integral (\ref{eqGauss})  nowadays is of a basic knowledge of Physics and  Mathematics (e.g. Electromagnetism, Fluid Dynamics).  
The integral is also expressed by one arrow among two circles, so-called (two-component) \emph{Gauss diagram formula}  
\[
\left\langle \LK~, G_L \right\rangle \quad.    
\]
It is the simplest two-component Gauss diagram formula.  In this paper, we study the second simple two-component Gauss diagram presenting  topological functions  
\[
\Ga, \quad \Gb \quad.
\]  
For these integer-valued functions, we call them \emph{multiple linking numbers}.  
The former is invariant of ambient  isotopy and further link homotopy of links.   Surprisingly, the other becomes  
a function showing the necessity of the second Reidemeister moves  (Theorem~\ref{main}~(\ref{ind})).  
\"{O}stlund \cite{Ostlund2001PhD} show that any one-component 
Gauss diagram formula 
cannot detect the necessity (Fact~\ref{factKnot}).  
In order to  recall the precise general statement (Fact~\ref{factKnot}) of \"{O}stlund, we pose Definition~\ref{KDdef}  here.  
\begin{definition}[knot/link diagram invariant (the knot case is given by Ostlund {\cite{Ostlund2001PhD}})]\label{KDdef}
We shall call a function of knot/link diagrams that is unchanged planar isotopy, but not necessarily by Reidemeister moves, a knot/link diagram invariant.   
\end{definition}
Since the $n$th Reidemeister move is often called $\Omega_n$-moves ($n=1, 2, 3$), so do we if no confusion occurs.    
\begin{fact}[{\cite[Chapter~IV, Theorem~2]{Ostlund2001PhD}}] \label{factKnot}
Let $v$ be a knot diagram invariant that is unchanged under $\Omega_1$- and $\Omega_3$-moves, and of finite degree (in $\Omega_2$).  Then $v$ is a knot invariant.  
\end{fact}
For links, this is not the case \cite{Ito2022_OnTypeII}.  Theorem~\ref{main}~(\ref{ind}) gives a simpler example  than that of \cite{Ito2022_OnTypeII}; moreover, in a case, the function in Theorem~\ref{main}~(\ref{ind}) counts the minimum number of $\Omega_2$-moves decreasing crossings (Corollary~\ref{CorSurjective}).  
\section{Main Result}\label{sec:main}
\begin{theorem}\label{main}
{\color{black}{
Let \[S =
\begin{picture}(100,0)
\put(0,-7){$\Ga$\quad,}
\end{picture}
T = 
\begin{picture}(100,0)
\put(0,-7){$\Gb$\quad. }
\end{picture}
\] Suppose that $L$ is an arbitrary  two-component link and $G_L$ is a Gauss diagram of $L$.   Then the following statements hold.   
\begin{enumerate}
\item The integer-valued function $\langle S, G_L \rangle$ is an isotopy (more precisely, link homotopy) invariant of $L$.   \label{inv}  
\item The integer-valued function $\langle T, G_L \rangle$ is invariant under $\Omega_1$- and $\Omega_3$-moves, whereas it changes under a $\Omega_2$-move.    \label{ind}
\end{enumerate}
}}
\end{theorem}
\begin{corollary}\label{CorSurjective}
For any nonpositive integer $n$, there exists a link diagram $D$ such that 
$\langle T, G_{D} \rangle = n$.   
\end{corollary}
If a $\Omega_2$-move decreases crossings, we call it a \emph{negative} $\Omega_2$-move.  
\begin{corollary}\label{CorDistance}
For any nonpositive integer $n$,   there exist two diagrams $D$ and $D'$ of $L$ such that $|\langle T , D' \rangle $ $-$ $\langle T , D \rangle|$ equals the minimum number of negative $\Omega_2$-moves among any  sequences of $\Omega_*$-moves $(*=1, 2, 3)$ from $D$ to $D'$.   
\end{corollary}
\begin{remark}
Seeing the proof in Section~\ref{sec:proof}, it is elementary to see that the extension of Theorem~\ref{main} to virtual link is straightforward.   
\end{remark}
\section{\color{black}{Preliminary}}\label{sec:Prelim}
\begin{figure}[htbp] 
   \centering
   \includegraphics[width=7cm]{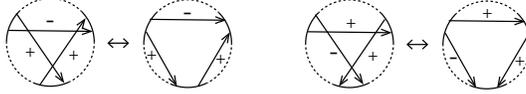} 
   \caption{The third Reidemeister move $\Omega_{3a}$ is presented by two moves of Gauss diagrams: $\Omega_{3+-++}$   (left) and $\Omega_{3+-+-}$ (right)}
   \label{OSTNito}
\end{figure}

\begin{definition}[cf.~\cite{ItoTakimura2020}]\label{rii}
Let $D$ be a link diagram.  Then Let $\rii(D)$ $=$  
the minimum number of $\Omega_2$-moves decreasing crossings in a sequence of $\Omega_*$-moves $(*=1, 3)$ and $\Omega_2$-moves decreasing crossings  from $D$ to $D'$.  
\end{definition}
In this paper, we freely use Polyak-Viro-\"{O}stlund's terminologies and notations \cite{Ostlund2004}, e.g. Gauss diagrams, arrow diagrams, Gauss or arrow diagram fragments, and Reidemeister moves \cite[Section~1.6, Section~4, and Table~1]{Ostlund2004} except for the replacement of $\Omega_{+---}$ in \cite[Table~1]{Ostlund2004} with $\Omega_{+-+-}$ in Figure~\ref{OSTNito}.  

By the definition of Gauss diagrams, we note the ambiguity representing a Gauss diagram when we do not fix the order of components and base points.  In our case, 

\begin{equation}\label{Represent}
\begin{split}
\begin{picture}(300,0)
\put(0,-17){$\Ga$} 
\put(75,-6){$\equiv$} 
\put(100,-17){$\Gal$} 
\put(175,-6){$\equiv$} 
\put(200,-17){$\GalS$} 
\put(275,-6){$\equiv$} 
\put(300,-17){$\GaS$, } 
\end{picture}\\
\\
\begin{picture}(300,0)
\put(0,-17){$\Gb$} 
\put(75,-6){$\equiv$} 
\put(100,-17){$\GbS$} 
\put(175,-6){$\equiv$} 
\put(200,-17){$\Gbl$} 
\put(275,-6){$\equiv$} 
\put(300,-17){$\GblS$} 
\end{picture}
\end{split}
\end{equation}
for 
\[\]
\[S =
\begin{picture}(100,0)
\put(0,-7){$\Ga$\quad,}
\end{picture}
T = 
\begin{picture}(100,0)
\put(0,-7){$\Gb$\quad.}
\end{picture}
\]
\begin{example}\label{eg_toruslink}
For ($2, 2n$)-torus standard link diagram with $2n$-positive crossing $G_{D(n)}$, 
\[
\langle S, G_{D(n)} \rangle = n^2, \langle T, G_{D(n)} \rangle = n(n-1).  
\]
Let $D'(n)$ be by applying exactly one $\Omega_2$-move to $D_n$ increasing crossings, 
\[
\langle S, G_{D'(n)} \rangle = n^2, \langle T, G_{D'(n)} \rangle = n(n-1)-1.  
\]
\end{example}

{\color{black}{
\section{\color{black}{Proof of Theorem~\ref{main}}}\label{sec:proof}
In this section, we give a proof of Theorem~\ref{main}.  We should  note that proofs of (\ref{inv}) and (\ref{ind}) of Theorem~\ref{main} are quite different.  The link invariance is given by checking a generating set of $\Omega_n$-moves ($n=1, 2, 3$), whereas we should see every pattern of $\Omega_3$-move to prove Theorem~\ref{main}~(\ref{ind}) because we cannot use $\Omega_2$-move.        

\subsection{Proof for (\ref{inv})}\label{proof_link}
Recall that $\{ \Omega_{1a}, \Omega_{1b}, \Omega_{2a}, \Omega_{3a} \}$ is a generating set \cite[Polyak (2010)]{Polyak2010}.   Note that it is clear that $\langle S, \cdot \rangle$ is invariant under Reidemeister moves in one component.  Note also that it is elementary to show that for any $\Omega_2$-move, $\langle S, \cdot \rangle$ does not change.  Thus, we focus on the $\Omega_{3a}$-move related to two components, which are represented by $\Omega_{+-+*}$ ($*=b, m, t$) of \cite[Table~1]{Ostlund2004} via Gauss diagrams.   

Each of $\Omega_{+-+*}$ ($*=b, m, t$) is presented by Gauss diagrams as 
\[ 
\begin{picture}(100,30)
\put(0,-7){$\OR$} 
\put(75,10){$\leftrightarrow$} 
\put(100,-7){$\ORS$}
\end{picture}
\]
up to signs 
and arrow-orientations.   By the same manner as the definition  \cite[Definition~5]{Ito2019} for curves or \cite[Definition~8]{ItoKotoriiTakamura2022} for virtual knots of \emph{relators}, every  relator corresponding to $\Omega_{+-+*}$ ($*=b, m, t$) is of type
\begin{equation}\label{relatorR}
\begin{split}
\left(
\begin{picture}(280,30)
\put(0,-17){$\OR$} 
\put(65,0){$+$} 
\put(80,-17){$\ORab$}
\put(145,0){$+$} 
\put(150,-17){$\ORar$}
\put(215,0){$+$} 
\put(230,-17){$\ORbr$}
\end{picture}
\right)
\\
- \left(
\begin{picture}(280,30)
\put(0,-17){$\ORS$} 
\put(65,0){$+$} 
\put(80,-17){$\ORabS$}
\put(145,0){$+$} 
\put(150,-17){$\ORarS$}
\put(215,0){$+$} 
\put(230,-17){$\ORbrS$}
\end{picture}
\right) \quad.  
\end{split}
\end{equation}
Let $R$ be the above type of relators.  Noting that the arrow diagram $S$ consists of two arrows, the argument returns to the fourth and fifth terms in $R$ corresponding to $\Omega_{+-+*}$ ($*=b, m, t$).  Then we find a three types subtractions, each of which consists of the fourth and fifth terms as  
\[  
\begin{picture}(100,30)
\put(-30, 10){$r_1 =$}
\put(0,-7){$\Aright$} 
\put(75,10){$-$} 
\put(100,-7){$\ArightS$, }
\end{picture}
\]
\[ 
\begin{picture}(100,30)
\put(-30, 10){$r_2 =$}
\put(0,-7){$\Bright$} 
\put(75,10){$-$} 
\put(100,-7){$\BrightS$,}  
\end{picture}
\] 
or
\[ 
\begin{picture}(100,30)
\put(-30, 10){$r_3 =$}
\put(0,-7){$\Aleft$} 
\put(75,10){$-$} 
\put(100,-7){$\AleftS$.}  
\end{picture}
\] 
In any case, $\langle S, r_i \rangle=0$ ($i=1, 2, 3$) noting    (\ref{Represent}) for $r_2$ (It is clear for either $r_1$ or $r_3$).  
The fact shows that these pairs induce invariances which we have been desired.  
$\hfill \Box$
\begin{remark}
In \cite[Proof of Theorem~2, case $\epsilon_4=1$]{Ito2019}, if you pose on $``|"$ on $\tau(z_0)$ to give a Gauss phrase, everything works to show this invariance by copying and pasting sentences; in this case, $R$ corresponds to \cite[Type SRI\!I\!I relator]{Ito2019}.  Thus we omit to do it.     
\end{remark}
\subsection{Proof of (\ref{ind})}
\label{proof_necessity}
Since it is clear that $\langle T, \cdot \rangle$ is invariant under $\Omega_1$-move, we shall check the difference values before/after applying a $\Omega_{3}$-move; specifically, we focus on $\Omega_{3a}$, $\Omega_{3b},\Omega_{3c}$, $\Omega_{3d}$, $\Omega_{3e}$, $\Omega_{3f}$, $\Omega_{3g}$, and $\Omega_{3h}$ by notations in \cite{Polyak2010}.  

Firstly, we consider $\Omega_{3a}$.  
Recalling (\ref{relatorR}) in Section~\ref{proof_link}, we focus on the subtraction $r_i$ ($i=1, 2, 3$ of  Section~\ref{proof_link}) that is the fourth minus fifth terms.  
Then it is elementary to see that $\langle T, r_i \rangle=0$, which implies the invariance under $\Omega_{3a}$.      
}}

Secondly, by the fact together with \cite[Table~2]{Ito2022_OnTypeII}, only the invariance of $\Omega_{3a}$ implies those of the other $\Omega_{3*}$ ($*=b, c, d, e, f, g, h$).  According to \cite{Polyak2010},  
every $\Omega_{3*}$ is generated by $\Omega_{3a}$ and $\Omega_{2*}$ as in Table~\ref{Table2} \cite[Table~2]{Ito2022_OnTypeII}.   As for the first case in Table~\ref{Table2}, an application of $\Omega_{3b}$ to $D_0$ to obtain a diagram $D_3$ is decomposed into the sequence 
\[
D_0 \stackrel{\Omega_{2c}}{\to} D_1  \stackrel{\Omega_{3a}}{\to} D_2  
\stackrel{\Omega_{2d}}{\to} D_3
\]
where $\Omega_{2c}$ ($\Omega_{3a}$, $\Omega_{2d}$,~resp.) is applied to $D_0$ ($D_1$, $D_2$,~resp.).  
Suppose that $\langle T, \cdot \rangle$ increasing $\alpha$ ($\beta$,~resp.) by this $\Omega_{2c}$ ($\Omega_{2d}$,~resp.).    Since $\langle T , G_{D_1} \rangle$ $=$ $\langle T , G_{D_2} \rangle$ by the invariance for $\Omega_{3a}$, then we follow differences step by step  
\[
\langle T , G_{D_0} \rangle \stackrel{+\alpha}{\to} \langle T , G_{D_1} \rangle  \stackrel{0}{\to} \langle T , G_{D_2} \rangle  
\stackrel{+\beta}{\to} \langle T , G_{D_3} \rangle.   
\]
It is clear that $(\alpha, \beta)=(0, 0)$ or $(-1, 1)$ which implies $\langle T , G_{D_0} \rangle$ $=$ $\langle T , G_{D_3} \rangle$.  Each of the other cases is similar to the above case of $\Omega_{3a}$.  

\begin{figure}[htbp] 
    \includegraphics[width=5cm]{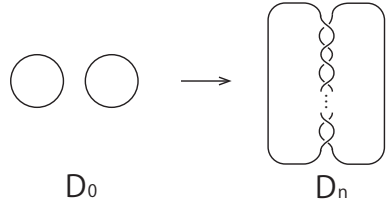} 
    \caption{The unlinked unknotted diagram $D_n$ with $2n$ crossings by applying $\Omega_2$-moves exactly $n$ times to two-component simple closed curves on the plane}
    \label{DiagN}
 \end{figure}

Finally, for a link diagram $D_n$ as in Figure~\ref{DiagN}, $\langle T, L_n \rangle$ $=$ ${}_n C_{2}$ $+$ ${}_n C_{2}$ $-n^2$ $=$ $-n$, which implies the necessity of $\Omega_2$-moves and also implies possibilities of changes of values.    
$\hfill \Box$
  \begin {table}
  \caption {$\Omega_{3*}$ returns to $\Omega_{3a}$
Each positive (negative,~resp.) $\Omega_{2*}$-move increases (decreases) crossings.  
  }
  \label {Table2}
  \begin {tabular}{lllll} \hline
    Type III 			& 	number of move	&	Positive Type II			& 	Type III			&	Negative Type II 		\\ \hline
    $\Omega_{3b}$	&		1		 	&	$\Omega_{2c}$  		&  	$\Omega_{3a}$		& 	$\Omega_{2d}$ 		\\ \hline
    $\Omega_{3c}$ 	& 		1			&	$\Omega_{2c}$  		& 	$\Omega_{3a}$	 	& 	$\Omega_{2d}$			\\ \hline
    $\Omega_{3d}$	& 		1			&	$\Omega_{2a}$ 		& 	$\Omega_{3b}$  	& 	$\Omega_{2b}$			\\
    $ $			&		2		 	&	$\Omega_{2c}$  		&  	$\Omega_{3a}$		& 	$\Omega_{2d}$			\\  \hline
    $\Omega_{3e}$	& 		1			& 	$\Omega_{2a}$ 		& 	$\Omega_{3b}$  	& 	$\Omega_{2b}$			\\ 
    $ $			&		2		 	&	$\Omega_{2c}$  		&  	$\Omega_{3a}$		& 	$\Omega_{2d}$			\\  \hline
    $\Omega_{3f}$ 	&  		1			&	$\Omega_{2d}$ 		&	$\Omega_{3a}$ 	& 	 $\Omega_{2c}$		\\  \hline
    $\Omega_{3g}$	& 		1			&	$\Omega_{2c}$ 		&	$\Omega_{3f}$   	&  	$\Omega_{2d}$			\\ 
    $ $			&		2			&	$\Omega_{2d}$ 		&	$\Omega_{3a}$		& 	 $\Omega_{2c}$		\\  \hline
    $\Omega_{3h}$ 	& 		1			&	$\Omega_{2a}$ 		&	$\Omega_{3f}$    	&  	$\Omega_{2b}$			\\
    $ $			&		2			&	$\Omega_{2d}$ 		&	$\Omega_{3a}$		& 	 $\Omega_{2c}$		\\  \hline
  \end {tabular}
  \end {table}

\noindent{\bf Proof of Corollary~\ref{CorSurjective}.}
\begin{proof}
By the computation in the last sentence of the  proof of Section~\ref{proof_necessity}, it is clear.  
\end{proof}

\noindent{\bf Proof of Corollary~\ref{CorDistance}.}
\begin{proof}
Let $D_0$ and $D_n$ be as in Figure~\ref{DiagN}.  It is sufficient to prove Claim~\ref{ClaimDn}.  
\begin{claim}\label{ClaimDn}
$D_n$ needs at least $n$ times $\Omega_2$-moves to obtain $D_0$, i.e. \[
|\langle T, D_n \rangle| = \rii(D_n). \]
\end{claim}
Assume that $\rii(D_{n}) \le n-1$.  
Note that any $\Omega_1$- and $\Omega_3$-moves do not affect increments/decrements of arrows between the different circles in the Gauss diagram $G_{D_n}$.  The $k$ times $\Omega_2$-moves cannot decrease values $|\langle T, \cdot \rangle|$ by at most $n-1$; that is, the $k$ times $\Omega_2$-moves cannot obtain $D_0$, which be in  contradiction to $\rii(D_{n}) \le n-1$.  Hence
\[
 \rii(D_{n}) \geq n.  
\]
By the construction of $D_{n}$, $\rii(D_{n}) \le n$.  Therefore, $
\rii(D_{n}) = n = |\langle T, D_{n} \rangle|$.  
\end{proof}
\section{Examples}
The integer-valued function derived from $S$ and $T$ can be applied into classical and virtual links. 

\subsection{Application of $S$ on a set of Virtual Links }\label{S_virtual}
\begin{figure}[htbp]
   \includegraphics[width=5in]{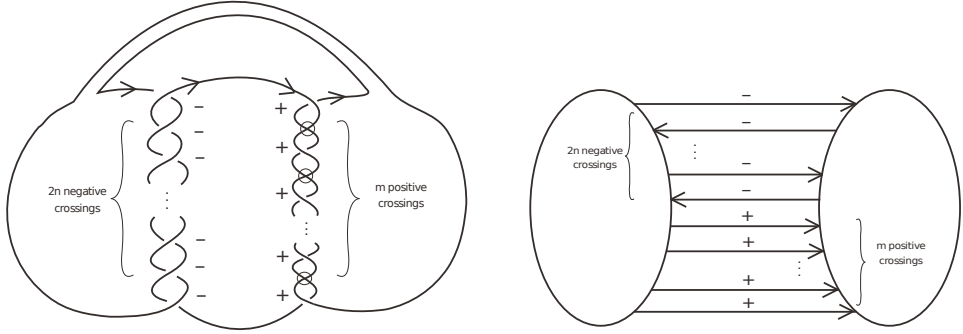} 
   \caption{Virtual link diagrams $\mathcal{L}_{m, n}$ and their Gauss diagrams $G_{\mathcal{L}_{m, n}}$}
   \label{egSujective}
\end{figure}

\begin{proposition}
Let $m \in \mathbb{Z}_{\geq 0}$ and $n \in \mathbb{Z}_{>0}$.  Let $\mathcal{L}_{m, n}$ be the set of $2$-component virtual links with $m$ virtual crossings and $2n$ positive crossings as in Figure~\ref{egSujective}.  Then if $m = n-1$, $\langle S, \mathcal{L}_{n-1, n} \rangle$  $\to$ $\mathbb{Z}_{> 0}$ is surjective.
\end{proposition}

\begin{proof}
Let $L_{m, n} \in \mathcal{L}_{m, n}$ and 
the Gauss diagram of $L_{m, n}$ is  as in Figure~\ref{egSujective}.  As we observed in a paragraph proceeding of Example~\ref{eg_toruslink}, we do not consider base point. Then
\begin{align*}
\langle S, L_{m, n} \rangle  &=  \left\langle \GaMinusMinus+ \GaSMinusMinus + \GaSMinusPlus , L_{m, n} \right\rangle \\
&= \frac{(n+1)n}{2}+\frac{n(n-1)}{2} -mn \\
&= n^2-mn = n(n-m).
\end{align*}
In particular, 
$\langle S, L_{n-1, n} \rangle$  = $n(n-(n-1))$ = $n$. 

Then the function $\langle S, \cdot \rangle$ from the set $\mathcal{L}_{n-1, n}$ of $2$-component virtual links to $\mathbb{Z}$ is surjective.
\end{proof}

\subsection{Application of $T$ on a set of Classical Links }\label{T_virtual}
\begin{figure}[htbp]
   \includegraphics[width=5in]{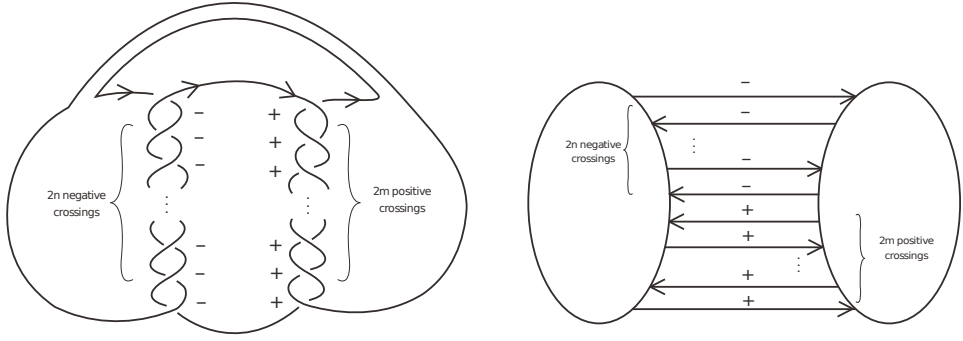} 
   \caption{Virtual link diagrams $\mathcal{K}_{m, n}$ and their Gauss diagrams $G_{\mathcal{K}_{m, n}}$}
   \label{egTsurjective}
\end{figure}

\begin{proposition}
For $n,m \in  \mathbb{Z}_{\geq 0}$, let $\mathcal{K}_{m,n}$ be the set of 2-component virtual links with 2n positive crossings and 2m negative crossings as in Figure~\ref{egTsurjective}.  Let $2\mathbb{Z}_{\leq 0}$ $=$ $\{ 2n~|~n \in \mathbb{Z}_{\leq 0} \}$.  
Then 
if m = n or n+1, $\langle T, \mathcal{K}_{m,n} \rangle$ $\to$ $2\mathbb{Z}_{\leq 0}$ is surjective.  
\end{proposition}

\begin{proof}
Let $K_{m,n} \in \mathcal{K}_{m,n}$ and the Gauss diagram is as in Figure~\ref{egTsurjective}.  As we observed in a paragraph proceeding of Example~\ref{eg_toruslink}, we do not consider base point. Then, 
\begin{align*}
\langle T, K_{m,n} \rangle  &= \Big \langle \GbMinusMinus + \GblMinusMinus + \GblPlusPlus + \GbPlusPlus \\
& +\GbMinusPlus + \GblMinusPlus, K_{m,n} \Big \rangle \\
&= \frac{n(n-1)}{2}+\frac{n(n-1)}{2} + \frac{m(m-1)}{2}+\frac{m(m-1)}{2} -mn -mn \\
&=n^2-n + m^2-m - 2mn \\
&=(n^2-2mn+m^2)-(n+m) \\
&= (n-m)^2-(n+m).  \\
\end{align*}

In particular,  $\langle T, K_{n, n} \rangle$ = $(n-n)^2-(n+n)$ = $-2n$. 
\\ \indent $\langle T, K_{n+1, n} \rangle$ = $(n-(n+1))^2-(n+(n+1))$ = $1-2n-1$ = $-2n$. 
\\ \indent Then the function $\langle T, \cdot \rangle$ from the set $\mathcal{K}_{n, n}$ and $\mathcal{K}_{n+1, n}$ of $2$-component virtual links to $2\mathbb{Z}_{\leq 0}$ is surjective.  
\end{proof}
 Note that links in the set $\mathcal{K}_{m, n}$ can obviously be resolved into the trivial link using $\Omega_2$-move starting from the center (the boundary of $2n$ and $2m$ crossings), but the function $\langle T, \cdot \rangle$ took a non-zero value. This indicates that the function  $\langle T, G_L \rangle$ is invariant under $\Omega_1$- and $\Omega_3$-moves, whereas it changes under $\Omega_2$-moves.
 


\bibliographystyle{plain}
\bibliography{RefTwo}
\end{sloppypar}
\end{document}